\documentclass{amsart}
\usepackage{amssymb}
\usepackage{amsmath}
\usepackage{amsfonts}

\theoremstyle{plain}

\numberwithin{equation}{section}
\input{tcilatex}

\begin{document}
\title[Semihyperrings Characterized by Their Hyperideals]{Semihyperrings
Characterized by Their Hyperideals}
\author{M. Shabir}
\address[M. Shabir]{Department of Mathematics Quaid-I-Azam University
Islamabad, Pakistan.}
\email[M. Shabir]{mshabirbhatti@yahoo.co.uk}
\author{Nayyar Mehmood}
\address[Nayyar Mehmood]{Department of Mathematics, School of Chemical and
Material Engineering (SCME), NUST, H-12, Islamabad, Pakistan.}
\email[Nayyar Mehmood]{nayyarnh@yahoo.co.uk}
\author{Piergiulio Corsini}
\address[P. Corsini]{Dept. of Biology and Agro-Industrial Economy, VIA DELLE
Scienze 208, 33100 UDINE, ITALY.}
\email{corsini2002@yahoo.com}
\urladdr{http://ijpam.uniud.it/journal/curriculum\_corsini.htm}
\date{March 13, 2010.}

\begin{abstract}
The concept of hypergroup is generalization of group, first was introduced
by Marty [9]. This theory had applications to several domains. Marty had
applied them to groups, algebraic functions and rational functions. M.
Krasner has studied the notion of hyperring in [11]. G.G Massouros and C.G
Massouros defined hyperringoids and apply them in generalization of rings in
[10]. They also defined fortified hypergroups as a generalization of
divisibility in algebraic structures and use them in Automata and Language
theory. T. Vougiouklis has defined the representations and fundamental
relations in hyperrings in [14, 15]. R. Ameri, H. Hedayati defined
k-hyperideals in semihyperrings in [2]. B. Davvaz has defined some relations
in hyperrings and prove Isomorphism theorems in [7]. The aim of this article
is to initiate the study of semihyperrings and characterize it with
hyperideals. In this article we defined semihyperrings, hyperideals, prime,
semiprime, irreducible, strongly irreducible hyperideals, m-systems,
p-systems, i-systems and regular semihyperrings. We also shown that the
lattice of irreducible hyperideals of semihyperring $R$ admits the structure
of topology, which is called irreducible spectrum topology.
\end{abstract}

\maketitle

\section{\protect\bigskip Introduction and Preliminaries}

The hyperstructures{\LARGE \ }are algebraic structures equipped with, at
least, one multivalued operation, called hyperoperation. The largest classes
of the hyperstructures are the ones called $Hv$\ -- structures. A\
hyperoperation *\ on a non-empty set $H$\ is a mapping from $H\times H$\ to
power set $P^{\ast }(H),$\ where $P^{\ast }(H)$\ denotes the set of all
non-empty subsets of $H$, that is $\cdot :$\ $H\times H$\ $\rightarrow
P^{\ast }(H):(x,y)\rightarrow x\cdot y\subseteq H$.

We denote here the hyperoperation as '$\cdot $'. If $A$ and $B$ are two
non-empty subsets of $H$, then

\qquad \qquad \qquad \qquad \qquad \qquad \qquad \qquad \qquad

\qquad \qquad \qquad \qquad \qquad \qquad $A\cdot B=\underset{\underset{%
{\large b\in B}}{a\in A}}{\cup }$ $a\cdot b$ \ \ \ 

\ \ \ \ if $A=\{a\}$ then$\qquad \qquad A\cdot B=\underset{b\in B}{\cup }%
a\cdot b$

\subsection{Definition:[2]}

A hyperoperation $\cdot $ on a non-empty set $H$ is called associative if $%
(a\cdot b)\cdot c=a\cdot (b\cdot c)$ for all $a,$ $b$ $,c\in H$.

\subsection{Definition:[2]}

The set $H$ with hyperoperation $\cdot $ is called hypergroupoid, we denote
it by $(H,\cdot ).$

\subsection{Definition:[2]}

A hypergroupoid $(H,\cdot )$ with associative hyperoperation is called a $%
semihypergroup.$

\bigskip

\subsection{Definition:[2]}

An element $e\in H$ is called an identity element of a semihypergroup $%
(H,\cdot )$ if for all $a\in H,$ we have\qquad $e\cdot a=a\cdot e=\{a\}.$

Sometime we can also write it, $\ e\cdot a=a\cdot e=a$ \ just for the sake
of easiness.

\subsection{Definition:[1]}

A non-empty set $H$ is called a canonical hypergroup or simply a hypergroup
if the hypergroupoid $(H,\cdot )$ satisfies the following properties

\begin{enumerate}
\item $(a\cdot b)\cdot c=a\cdot (b\cdot c)$ for all $a,$ $b$ $,c\in H$.

\item An element $e$ exists in $H$, such that $e\cdot a=a\cdot e=a$ for all $%
a\in H.$ That is $e$ is an identity element of $H.$

\item For all $x\in H,$ there exists a unique $x^{-1}\in H,$such that $e\in
x\cdot x^{-1}=x^{-1}\cdot x.$

\item $(H,\cdot )$ is reversible, that is if $z\in x\cdot y,$ then $x\in
z\cdot y^{-1}$ and $y\in x^{-1}\cdot z.$
\end{enumerate}

\subsection{Definition:[2]}

A hypergroup is a semihypergroup $(H,\cdot )$\ on which the reproduction
axiom is valid, that is $a\cdot H=H\cdot a=H,$\ for all $a$ in $H$.

\subsection{Definition:[1]}

A hyperring is an algebraic structure $(R,+,\cdot )$ which satisfies the
following axioms:

\begin{enumerate}
\item $(R,+)$ is an abelian canonical hypergroup;

\item $(R,\cdot )$ is a semigroup having $0$ as an absorbing element, that
is, $0\cdot a=a\cdot 0=0$ for all $a\in R.$

\item The multiplication is distributive with respect to hyperoperation +.
\end{enumerate}

A canonical hypergroup means the associativity and commutativity of the
hyperaddition, the existence of a unique $-a$ for $a\in R,$ such that $0\in
a+(-a)$ and the implication: $a\in b+c$ implies $b\in a+(-c).$

The following elementary facts about canonical hypergroups follow easily
from the above axioms:

\begin{enumerate}
\item $-(-a)=a$, for all $a\in R;$

\item $a+0=a,$ for all $a\in R;$

\item $0$ is unique such that for every $a\in R,$ there is an element $(-a)$
with property $0\in a+(-a);$

\item $a+R=R$, for all $a\in R;$

\item $-(a+b)=-a-b,$ or all $a,$ $b\in R.$
\end{enumerate}

\subsection{Example:[1]}

Let $(A,+,\cdot )$ be a ring and $N$ be a normal subgroup of its
multiplicative semigroup. Then the multiplicative classes $\bar{x}=xN(x\in
A) $ form a partition of $R,$ and let $\bar{A}=A/N$ be the set of these
classes. The product of $\bar{x},$ $\bar{y}\in \bar{A}$ as subsets of $A$ is
again a class $(\func{mod}N),$ and their sum as such subset is a union of
such classes. If we define the product $\bar{x}o\bar{y}$ in $\bar{A}$ of $%
\bar{x},$ $\bar{y}$ in $\bar{A}$ as equal to their product as subset of $A,$
and their sum $\bar{x}\oplus \bar{y}$ in $\bar{A}$ as the set of all $\bar{z}%
\in \bar{A}$ contained in their sum as subsets of $A,$ that is $\bar{x}%
\oplus \bar{y}=\{\bar{z}$ $%
{\vert}%
$ $z\in \bar{x}+\bar{y}\},\qquad $and $\bar{x}o\bar{y}=\overline{x\cdot y}$,
then $(A,\oplus ,\cdot )$ is a hyperring.

\section{HYPERIDEALS IN SEMIHYPERRING}

\subsection{Definition:}

A semihyperring is an algebraic structure $(R,+,\cdot )$ which satisfies the
following axioms:

$(1)$ $(R,+)$ is a commutative hypermonoid, that is;

$\qquad (i)$ $(x+y)+z=x+(y+z)$ for all $x,y,z\in R$

$\qquad (ii)$ There is $0\in R,$ such that\ $x+0=0+x=0$\ for all $x$ $\in R$

$\qquad (iii)$ $x+y=y+x\ for$ all $x,y$ $\in $ $R.$

$(2)$ $(R,\cdot )$ is semigroup, that is$\ (i)$ $\ \ x\cdot (y\cdot
z)=(x\cdot y)\cdot z$ \ for all $x,y,z\in R$

$(3)$ The multiplication is is distributive with respect to hyperoperation
'+' that is;

$\qquad x\cdot (y+z)=x\cdot y+x\cdot z$ and $(x+y)\cdot z=x\cdot z+y\cdot z$
for all $x,$ $y,$ $z\in R.$

$(4)$ The element $0$ $\in $ $R$, is an absorbing element, that is;

$\qquad \ x\cdot 0=0\cdot x=0$\ for all\ \ $x\in R$

\subsection{Example:}

Let $X$ be a nonempty finite set and $\tau $ is a topology on $X.$ We define
the hyperoperations of the addition and the multiplication on $\tau $ as;

For any $A,$ $B\in \tau ,$ $A+B=\{A\cup B\}$ and $A\cdot B=A\cap B.$ Then $%
(\tau ,+,\cdot )$ is a semihyperring with absorbing element and additive

identity $\phi $ and multiplicative identity $X.$

\subsection{Example:}

Let us consider a set $S=\left\{ \left[ 
\begin{array}{cc}
a & b \\ 
c & d%
\end{array}%
\right] \text{ }%
{\vert}%
\text{\ }a,\text{ }b,\text{ }c,\text{ }d\in W\right\} ,$ where $W$ is a set
of whole numbers. We define the hyperoperation

of addition and multiplication as;

For $A=\left[ 
\begin{array}{cc}
a & b \\ 
c & d%
\end{array}%
\right] $ and $B=\left[ 
\begin{array}{cc}
a & b \\ 
c & d%
\end{array}%
\right] $ be taken from $S.$

$A+B=\left[ 
\begin{array}{cc}
a & b \\ 
c & d%
\end{array}%
\right] +\left[ 
\begin{array}{cc}
\acute{a} & \acute{b} \\ 
\acute{c} & \acute{d}%
\end{array}%
\right] =\left\{ \left[ 
\begin{array}{cc}
a+\acute{a} & b+\acute{b} \\ 
c+\acute{c} & d+\acute{d}%
\end{array}%
\right] \right\} \subseteq S$ and

$A\cdot B=%
\begin{array}{cc}
a & b \\ 
c & d%
\end{array}%
\cdot 
\begin{array}{cc}
\acute{a} & \acute{b} \\ 
\acute{c} & \acute{d}%
\end{array}%
$

Then $(S,+,\cdot )$ is a semihyperring with additive identity as a Null
matrix and multiplicative identity as a identity matrix.

\subsection{Definition:}

A semihyperring $(R,+,\cdot )$ is called commutative if \ $a\cdot b=b\cdot a$
for all $a,$ $b\in R.$

\subsection{Definition:}

A semihyperring $(R,+,\cdot )$ is called a semihyperring with identity if $%
1_{R}$ $\in R,$ such that $x\cdot 1_{R}=1_{R\text{ \ }}\cdot x=x$ for all $%
x\in R.$

\subsection{Definition:}

An element $x\in R$ is called unit if and only if there exists $y\in R$ such
that,\ $1_{R}=x\cdot y=y\cdot x.$ The set of all unit elements of a

semihyperring $R$ is denoted by $U(R).$

\subsection{Definition:}

A non-empty subset $S$ of a semihyperring $(R,+,\cdot )$ is called a
subsemihyperring of $R$ if

$(i)$ $a+b\subseteq S$ for all $a,$ $b\in S.$

$(ii)$ $a\cdot b\in S$ for all $a,$ $b\in S.$

\subsection{Definition:}

A left hyperideal of a semihyperring $R$ is a non-empty subset $I$ of $R$,
satisfying

(i) \qquad For all $x,$ $y\in I$, $x+y\subseteq I$

(ii)\qquad For all $a\in I$ and $x\in R,$ $x\cdot a\in I.$

A right hyperideal of a semihyperring $R$ is defined in a similar way.

\subsection{Definition:}

A hyperideal of a semihyperring $R$ is a non-empty subset of $R$, which is
both a left and a right hyperideal of $R.$

An absorbing element $0\in R$ also belongs to any hyperideal $I$ of $R.$
Since $0\in R$ so for any $a\in I$, we have

$0\cdot a=0\in I.$

\subsection{Proposition:}

A left hyperideal of a semihyperring is a subsemihyperring.

\begin{proof}
Let $I$ be a left hyperideal of a semihyperring $R.$

For any $a,$ $b\in I,$ we have $a+b\subseteq I.$

Also for any $a,$ $b\in I,$ since $a\in I\subseteq R.$ We have $a\cdot b\in
I.$

Hence $I$ is a subsemihyperring of the semihyperring $R.$
\end{proof}

\subsection{Definition:}

For every $x\in R$ there exists one and only one -$x\in R$ such that $0\in
x+($-$x)=($-$x)+x.$We shall write $\hat{x}$ for -$x$ and call it the
opposite of $x.$

Denote the set of all apposite elements of a semihyperring $R$, by $%
V_{H}(R), $ that is;

\qquad $V_{H}(R)=\{a\in R$ $%
{\vert}%
$ $\exists $ $b\in R,$ $and$ $o\in a+b\}$

\subsection{Definition:}

A nonempty subset $A$ of a semihyperring $R$ is called semihypersubtractive
if and only if $a\in A\cap V_{H}(R)$ implies that $\hat{a}\in A\cap
V_{H}(R), $ where $\hat{a}$ is the opposite of $\ a.$

\subsection{Definition:}

A nonempty subset $A$ of a semihyperring is called hypersubtractive if and
only if for $a\in A,$ and $a+b\subseteq A,$ implies $b\in A.$

\subsection{Definition:}

A nonempty subset $A$ of a semihyperring is called hyperstrong if and only
if $a+b\subseteq A,$ implies $a\in A$ and $b\in A.$

\subsection{Definition:}

A semihyperring $R$ is called additively reversive if it satisfy the
reversive property with respect to the hyperoperation of addition that is,
for $a\in b+c,$ implies $b\in a+\hat{c},$ and $c\in a+\hat{b}.$

\subsection{Theorem:}

A semihyperring is a hyperring if and only if $V_{H}(R)=R$ and $R$ is
additively reversive.

\begin{proof}
Suppose $(R,+,\cdot )$ is a hyperring. Clearly $V_{H}(R)\subseteq R\qquad
\qquad \rightarrow $ $(i)$

For any $a\in R,$ since $(R,+,\cdot )$ is a hyperring so there exists $\hat{a%
}\in R,$ such that \ $0\in a+\hat{a}.$ Thus $a\in V_{H}(R).$ So $R\subseteq
V_{H}(R)\qquad \qquad \qquad \qquad \rightarrow (ii)$

From $(i)$ and $(ii)$, we have $V_{H}(R)=R.$

Also as $R$ is a hyperring so is additively reversive.

$Conversely$ suppose that $R$ is additively reversive semihyperring and $%
V_{H}(R)=R.$

Let $s\in R=V_{H}(R),$ so there exists an element $\hat{s}\in V_{H}(R)=R,$
such that $0\in s+\hat{s}.$

Also since $R$ is additively reversive for $a\in b+c,$ implies $b\in a+\hat{c%
},$ and $c\in a+\hat{b}.$

Hence $(R,+)$ is a canonical hypergroup, so $R$ is a hyperring.
\end{proof}

\subsection{Definition:}

A semihyperring $(R,+,\cdot )$ is called zero sumfree if and only if $o\in r+%
\acute{r}$ implies $r=\acute{r}=o,$ or equivalently a semihyperring $%
(R,+,\cdot )$ is called zero sumfree if and only if $o\in r+\acute{r}$
implies $r\in \{o\},$ and $\acute{r}\in \{o\}.$

\subsection{Theorem:}

A semihyperring $(R,+,\cdot )$ is zero sumfree if and only if

$V_{H}(R)=\{o\}.$

\begin{proof}
Suppose $V_{H}(R)=\{o\}$

To prove $R$ is zero sumfree.

Let $a,$ $b\in R,$ such that $o\in a+b,$ or $o\in b+a,$ implies $a,$ $b\in
V_{H}(R)=\{o\},$ implies $a\in \{o\},$ and $b\in \{o\},$

that is $a=b=\{o\}.$ $R$ is zero sumfree.

$Conversely$ suppose that $R$ is zero sumfree. Since $o\in o+o,$ implies $%
\{o\}\subseteq V_{H}(R).$

Now let $r\in V_{H}(R),$ then there exist $\acute{r}\in R,$ such that $o\in
r+\acute{r}.$ But $R$ is zero sumfree, so $r=\acute{r}=o,$ implies $r\in
\{o\},$ and this implies that $V_{H}(R)\subseteq \{o\}.$

Hence $V_{H}(R)=\{o\}.$
\end{proof}

\subsection{Remark:}

Let $B$ be a hypersubtractive subset of a semihyperring $(R,+,\cdot )$, then
for any $b\in B,$ we have $o+b\subseteq B,$ since $B$ is hypersubtractive
so, $o+b\subseteq B,$ implies $o\in B.$

\subsection{Remark:}

Every hyperstrong subset of a semihyperring $R$ is hypersubtractive.

\begin{proof}
Suppose $A$ is a hyperstrong subset of a semihyperring $R$.

Let $a\in A,$ and $a+b\subseteq A.$ Since A is hyperstrong, so $b\in A.$

Hence $A$ is hypersubtractive.
\end{proof}

\subsection{Remark:}

A hypersubtractive subset $B$ of a semihyperring $R$ is semihypersubtractive
if for any $a\in B$ implies $a+\acute{a}=\{o\}.$

\begin{proof}
Let $B$ be a hypersubtractive subset of semihyperring $R,$ and for any $a\in
B$ implies $a+\acute{a}=\{o\}$

Let $a\in B\cap V_{H}(R),$ implies $a\in B,$ and $a\in V_{H}(R).$ Since $B$
is hypersubtractive so $o\in B,$ and for any $a\in B$ implies $a+\acute{a}%
=\{o\}.$ This implies that

\qquad \qquad \qquad \qquad $\{o\}=a+\acute{a}\subseteq B,$ implies $\acute{a%
}+a\subseteq B,$ and this implies, $\acute{a}\in B.$

Also $a\in V_{H}(R),$and we have $o\in \acute{a}+a,$ implies $\acute{a}\in
V_{H}(R).$

Hence $\acute{a}\in B\cap V_{H}(R),$ so $B$ is semihypersubtractive.
\end{proof}

\subsection{Theorem:}

Let $(R,+,\cdot )$ is a semihyperring and $\{I_{i}\}_{i\in \wedge }$ be a
family of hyperideals of $R.$ Then $\underset{_{i\in \wedge }}{\cap }I_{i}$
is also a hyperideal of $R.$

\begin{proof}
Let $a,$ $b\in $ $\underset{_{i\in \wedge }}{\cap }I_{i},$ then $a,$ $b\in
I_{i}$ for all $i\in \wedge .$ Since each $I_{i}$ is a hyperideal so $%
a+b\subseteq I_{i}$ for all $i\in \wedge .$ Thus $a+b\subseteq \underset{%
_{i\in \wedge }}{\text{ }\cap }I_{i}.$

Now for $x\in R$ , and $a\in $ $\underset{_{i\in \wedge }}{\cap }I_{i}.$
Since $a\in $ $\underset{_{i\in \wedge }}{\cap }I_{i},$ implies $a\in I_{i}$
for all $i\in \wedge ,$ and each $I_{i}$ is a hyperideal so $x\cdot a\in
I_{i}$ for all $i\in \wedge ,$

implies $x\cdot a\in \underset{_{i\in \wedge }}{\text{ }\cap }I_{i}.$ Again
for $x\in R$ , and $a\in $ $\underset{_{i\in \wedge }}{\cap }I_{i}.$ Since $%
a\in $ $\underset{_{i\in \wedge }}{\cap }I_{i},$ implies $a\in I_{i}$ for
all $i\in \wedge ,$ and each $I_{i}$ is a

hyperideal so $a\cdot x\in I_{i}$ for all $i\in \wedge ,$ implies $a\cdot
x\in \underset{_{i\in \wedge }}{\text{ }\cap }I_{i}.$

Hence $\underset{_{i\in \wedge }}{\cap }I_{i}$ is a hyperideal of a
semihyperring $R.$

Hence intersection of hyperideals of a semihyperring $R$ is hyperideal of $%
R. $
\end{proof}

\subsection{Theorem:}

If $A$ and $B$ are hyperideals of a semihyperring $R$. Then $A+B$ is the
smallest hyperideal of $R,$ containing both $A$ and $B$.

Where $A+B=\underset{%
\begin{array}{c}
a\in A \\ 
b\in B%
\end{array}%
}{\dbigcup }\left( a+b\right) .$

\begin{proof}
$(i)$ Let $x,$ $y\in A+B,$ then there exist $a_{1},$ $a_{2}\in A,$ and $%
b_{1},$ $b_{2}\in B,$

such that, $x\in a_{1}+b_{1}$ and $y\in a_{2}+b_{2}.$

Now\qquad $x+y\subseteq (a_{1}+b_{1})+(a_{2}+b_{2}),$ and by using
associativity and commutativity of $(R,+,\cdot ).$ We have$\qquad
x+y\subseteq (a_{1}+a_{2})+(b_{1}+b_{2}),$ since $A$ and $B$ are hyperideals
of $R,$ so $(a_{1}+a_{2})$ $\subseteq A$, $(b_{1}+b_{2})\subseteq B.$ This
implies $(a_{1}+a_{2})+(b_{1}+b_{2})\subseteq A+B.$

So $x+y\subseteq A+B.$

$(ii)$ Let $r\in R,$ and $x\in A+B,$ then there exist $a\in A,$ $b\in B,$
such that, $x\in a+b.$

Now consider,

\qquad \qquad \qquad $r\cdot x\in r\cdot (a+b)=r\cdot a+r\cdot b,$ by
distributivity of $R.$

Since $A$ and $B$ are hyperideals of semihyperring $R,$ so for any $a\in A,$ 
$b\in B$ and $r\in R,$ we have $r\cdot a\in A,$ $r\cdot b\subseteq B,$
implies $r\cdot a+r\cdot b\subseteq A+B,$ this implies that $r\cdot x\in
A+B. $

Again for any $r\in R,$ and $x\in A+B,$ then there exist $a\in A,$ $b\in B,$
such that, $x\in a+b.$

Also we have $x\cdot r\in (a+b)\cdot r=a\cdot r+b\cdot r,$ by distributivity
of $R$ since $A$ and $B$ are hyperideals of semihyperring $R.$

So for any $a\in A,$ $b\in B$ and $r\in R,$ we have $a\cdot r\in A$, $b\cdot
r\in B,$ implies $a\cdot r+b\cdot r\subseteq A+B,$

this implies that $x\cdot r\in A+B.$ Hence $A+B$ is the hyperideal of $R.$

Now we will show that $A+B$ is the smallest hyperideal of $R$ containing $%
A\cup B.$

First we will prove $A\cup B\subseteq A+B,$ for this let $x\in A\cup B,$
then either $x\in A$ or $x\in B.$

Since $A$ and $B$ are hyperideals of semihyperring $R.$ So $0\in A,$ and $%
0\in B,$ where $0$ is an absorbing element of semihyperring $R.$

Now if $x\in A,$ since $x=x+0\subseteq A+B,$ for $0\in B,$ implies $x\in
A+B. $ Also if $x\in B,$ for $0\in A,$ we have $x=0+x\subseteq A+B,$

so in both cases $x\in A+B.$ Hence $A\cup B\subseteq A+B.$

Now to prove $A+B$ is the smallest hyperideal, let $K$ be any other
hyperideal of $R,$ containing both $A,$ $B.$

To prove $A+B\subseteq K.$

Let $x\in A+B,$ then there exist $a\in A,$ $b\in B,$ such that, $x\in a+b.$
Since $a\in A\subseteq A\cup B\subseteq K,$ implies $a\in K,$ also as $b\in
B\subseteq A\cup B\subseteq K.$ Since $a,$ $b\in K,$ and $K$ is hyperideal
of $R,$ so $a+b\subseteq K.$ Hence $A+B\subseteq K.$

Hence $A+B$ is the smallest hyperideal of $R,$ containing both $A$ and $B.$
\end{proof}

\subsection{Lemma:}

If $R$ is a semihyperring with unity and $a\in R$ then

$a\cdot R=\underset{r\in R}{\cup }\left\{ \underset{finit}{\sum }a\cdot
r\right\} $ $\left( R\cdot a=\underset{r\in R}{\cup }\left\{ \underset{finit}%
{\sum }r\cdot a\right\} \right) $ is the smallest

right(left) hyperideal of $R$ containing $a.$

\begin{proof}
Let $x,y\in a\cdot R,$ then $x\in \underset{finit}{\sum }a\cdot r_{i},$ and $%
y\in \underset{finit}{\sum }a\cdot r_{j}.$

So we have $x+y\subseteq a\cdot \underset{finit}{\sum }(r_{i}+r_{j})%
\subseteq a\cdot R,$ implies $x+y\subseteq a\cdot R.$

Now for any $x\in a\cdot R,$ and $s\in R,$ consider $x\cdot s\in (\underset{%
finit}{\sum }a\cdot r)\cdot s=\underset{finit}{\sum }a\cdot (r\cdot
s)\subseteq \underset{r\in R}{\cup }\left\{ \underset{finit}{\sum }a\cdot
r\right\} =a\cdot R,$

implies $x\cdot s\in a\cdot R.$

So $a\cdot R$ is right hyperideal of semihyperring $R.$

Now to prove $a\cdot R$ is the smallest hyperideal containing $a.$ For this,
let $K$ be any other right hyperideal

of $R$ containing $a.$ Then,

\qquad \qquad \qquad \qquad $\underset{finit}{\sum }a\cdot r_{i}\subseteq K,$
for all $r_{i}\in R,$ implies $\underset{r\in R}{\cup }\left\{ \underset{%
finit}{\sum }a\cdot r\right\} \subseteq K,$ this

implies $a\cdot R\subseteq K.$ Also $1_{R}\in R,$ any $a=a\cdot 1_{R}\in 
\underset{r\in R}{\cup }\left\{ \underset{finit}{\sum }a\cdot r\right\}
=a\cdot R$,

implies $a\in a\cdot R.$ Hence $a\cdot R=\underset{r\in R}{\cup }\left\{ 
\underset{finit}{\sum }a\cdot r\right\} $ is the smallest right hyperideal
of $R$ containing $a\in R.$
\end{proof}

\bigskip

Let $X$ be nonempty subset of a semihyperring $R$. Then the hyperideal
generated by $X$ is the intersection of all hyperideals of $R$ which
contains $X$, and we denote it by $\left\langle X\right\rangle $. Hence we
have,

\qquad \qquad \qquad \qquad \qquad $\left\langle X\right\rangle =\cap \{I$ $%
| $ $X\subseteq I,and$ $I$ $is$ $hyperideal$ $of$ $R\}$

\subsection{Lemma:}

If $R$ is commutative semihyperring with unity, and $a\in R\backslash \{0\}.$
Then $\left\langle a\right\rangle =\underset{r\in R}{\cup }\left\{ \underset{%
finit}{\sum }a\cdot r\right\} .$

\subsection{Proposition:}

If $S$ is a subsemihyperring of a semihyperring $R$ and $I$ is a hyperideal
of $R$ then;

\begin{enumerate}
\item $S+I$ is a subsemihyperring of $R.$

\item $S\cap I$ is a hyperideal of $S.$
\end{enumerate}

\begin{proof}
$(1)\qquad S+I=\underset{\underset{i\in I}{s\in S}}{\cup }\left( s+i\right) $

Obviously for $0\in S,$ and $0\in I,$ we have $\{0\}=0+0\subseteq S+I.$

implies $0\in S+I$ so $S+I\neq \phi .$

Let $x,$ $y\in S+I,$ then there exist $s_{1},$ $s_{2}\in S,$ and $i_{1},$ $%
i_{2}\in I$ such that,

$x\in s_{1}+i_{1},$ and $y\in s_{2}+i_{2},$ then

$x+y\subseteq \left( s_{1}+i_{1}\right) +\left( s_{2}+i_{2}\right) =\left(
s_{1}+s_{2}\right) +\left( i_{1}+i_{2}\right) \subseteq S+I.$

Now $x\cdot y\in \left( s_{1}+i_{1}\right) \cdot \left( s_{2}+i_{2}\right)
=s_{1}s_{2}+s_{1}i_{2}+i_{1}s_{2}+i_{1}i_{2}$

\qquad \qquad \qquad \qquad \qquad \qquad \qquad $=s_{1}s_{2}+\left(
s_{1}i_{2}+i_{1}s_{2}+i_{1}i_{2}\right) \subseteq S+I,$ since $I$ is a
hyperideal of $R,$

so $s_{1}i_{2}+i_{1}s_{2}+i_{1}i_{2}\subseteq I.$

Hence $S+I$ is a subsemihyperring of $R.$

$(2)\qquad $Let $a,$ $b\in S\cap I,$ implies $a,$ $b\in S$ and $a,$ $b\in I.$

Since $S$ is a subsemihyperring and $I$ is a hyperideal so,

$a+b\subseteq S,$ and $a+b\subseteq I,$ which implies $a+b\subseteq S\cap I.$

Now let $s\in S,$ then obviously $s\cdot a\in S,$ also $s\cdot a\in I,$
since $I$ is hyperideal and $S$ is subsemihyperring.

Also for $s\in S,$ we have $a\cdot s\subseteq S,$ also $a\cdot s\subseteq I,$
since $I$ is hyperideal and $S$ is subsemihyperring.

Hence $H\cap I$ is a hyperideal of $S.$
\end{proof}

\subsection{Definition:}

Let $A$ be a nonempty subset of semihyperring $R$ with unity. Then the set

$(O:A)_{H}=\{r\in R%
{\vert}%
r\cdot a=o,$ for all $a\in A\}$ is called left annihilator hyperideal of $A.$

The right annihilator hyperideal is defined as;\qquad $(A:O)_{H}=$ $\{r\in R%
{\vert}%
a\cdot r=o,$ for all $a\in A\}$

\subsection{Proposition:}

$(O:A)_{H}=\{r\in R%
{\vert}%
r\cdot a=o,$ for all $a\in A\}$ is left hyperideal of $R.$

\begin{proof}
To prove $(O:A)_{H}=\{r\in R%
{\vert}%
r\cdot a=o,$ for all $a\in A\}$ is left annihilator hyperideal of $R.$

Obviously $o\in (O:A)_{H}$ since $o\cdot a=o,$ for all $a\in A,$ implies $%
(O:A)_{H}\neq \phi .$

Let $r_{1},$ $r_{2}\in (O:A)_{H}$ implies $r_{1}\cdot a=o$ and $r_{2}\cdot
a=o,$ for all $a\in A.$

Since $o=o+o=r_{1}\cdot a+r_{2}\cdot a=(r_{1}+r_{2})\cdot a,$ implies $%
r_{1}+r_{2}\subseteq (A:O)_{H}$ for all $a\in A.$

Now let $r_{1}\in (O:A)_{H},$ implies $r_{1}\cdot a=o$ for all $a\in A.$
Then for any $r\in R.$

Consider $(r\cdot r_{1})\cdot a=r\cdot (r_{1}\cdot a)=r\cdot o=o,$ since $%
r_{1}\cdot a=o$ for all $a\in A.$

So $r\cdot r_{1}\in $ $(O:A)_{H},$ for any $r\in R,$ and $r_{1}\in $ $%
(O:A)_{H}.$ Hence $(O:A)_{H}$ is left hyperideal of $R.$
\end{proof}

\subsection{Theorem:}

If $A$ is right (left) hyperideal of a semihyperring $R.$ Then $(A:O)_{H}$ $%
\left( (O:A)_{H}\right) $\ is a hyperideal of $R.$

\begin{proof}
Suppose $A$ is right hyperideal of a semihyperring $R,$ and

$(A:O)_{H}=\{r\in R%
{\vert}%
a\cdot r=o,$ for all $a\in A\}$

Let $r_{1},$ $r_{2}\in (A:O)_{H}$ implies $a\cdot r_{1}=o$ and $a\cdot
r_{2}=o,$ for all $a\in A.$

Since $o=o+o=a\cdot r_{1}+a\cdot r_{2}=a\cdot (r_{1}+r_{2}),$ for all $a\in
A $ implies $r_{1}+r_{2}\subseteq (A:O)_{H}$.

Now let $r_{1}\in (A:O)_{H},$ implies $a\cdot r_{1}=o$ for all $a\in A.$
Then for any $r\in R,$ consider $a\cdot (r\cdot r_{1})=(a\cdot r)\cdot r$

since $a\in A$ and $A$ is right hyperideal so $a\cdot r\in A,$ for all $r\in
R.$Since $r_{1}\in $ $(A:O)_{H},$ so we have

$a\cdot (r\cdot r_{1})=(a\cdot r)\cdot r_{1}=o,$ implies $r\cdot r_{1}\in
(A:O)_{H}.$

Now consider for any $r\in R,$ and $r_{1}\in $ $(O:A)_{H},$ $a\cdot
(r_{1}\cdot r)=(a\cdot r_{1})\cdot r=o.$ Since $a\cdot r_{1}=o$ for all $%
a\in A.$

So $r_{1}\cdot r\in (A:O)_{H},$ for any $r\in R,$ and $r_{1}\in $ $%
(A:O)_{H}. $

Hence $(A:O)_{H}$ is hyperideal of $R.$
\end{proof}

\subsection{Corollary:}

If $R$ is commutative semihyperring then $(O:A)_{H}$ and $(A:O)_{H}$ are
hyperideals.

\begin{proof}
Straight forward.
\end{proof}

\bigskip

\section{MULTIPLICATIVELY REGULAR SEMIHYPERRINGS}

\subsection{\protect\bigskip Definition:}

An element $a\in R$ is called multiplicatively regular if and only if there
exists an element $b\in R$ such that $\ a=a\cdot b\cdot a.$

A semihyperring $R$ is called multiplicatively regular if and only if each
element of $R$ is multiplicatively regular.

\subsection{Definition:}

If $I$ and $J$ are two hyperideals of semihyperring $R$. Then we define the
product $IJ$ of hyperideals $I$ and $J$ as;\qquad $IJ=\underset{b_{i}\in J}{%
\underset{a_{i}\in I}{\cup }}\left( \underset{finit}{\sum }a_{i}\cdot
b_{i}\right) $

\subsection{Proposition:}

The following conditions on a semihyperring $R$ with unity are equivalent.

\begin{enumerate}
\item $R$ is multiplicatively regular.

\item $H\cap I=HI$ for every left hyperideal $I$ and right hyperideal $H$ of
semihyperring $R.$
\end{enumerate}

\begin{proof}
$(1)\Longrightarrow (2)$

Suppose that, the semihyperring $R$ is multiplicatively regular.

To prove $H\cap I=HI$ for every left hyperideal $I$ and right hyperideal $H$
of semihyperring $R.$

Let $x\in HI=\underset{b_{i}\in I}{\underset{a_{i}\in H}{\cup }}\left( 
\underset{finit}{\sum }a_{i}\cdot b_{i}\right) ,$ implies $x\in \cup X_{i},$
where each

$X_{i}=\underset{finit}{\sum }a_{i}\cdot b_{i}$ for $a_{i}$ $\in H$ and $%
b_{i}\in I.$ So $x\in X_{j}$ for some $j.$

So $x\in X_{j}=\underset{finit}{\sum }a_{i}\cdot b_{i}$ for $a_{i}$ $\in H$
and $b_{i}\in I.$ Since for each $a_{i}\in H,$ we have $a_{i}\cdot b_{i}\in
H,$ for $b_{i}\in R,$ because $H$ is right hyperideal of semihyperring $R.$

Hence $\underset{finit}{\sum }a_{i}\cdot b_{i}\subseteq H,$ implies $x\in
X_{j}=\underset{finit}{\sum }a_{i}\cdot b_{i}\subseteq H,$ we have $x\in H,$
so $HI\subseteq H.$

Also for each $b_{i}\in I,$ we have $a_{i}\cdot b_{i}\in I,$ (for each $i$ )
for $a_{i}\in R,$ because $I$ is

left hyperideal of semihyperring $R.$ Hence $\underset{finit}{\sum }%
a_{i}\cdot b_{i}\subseteq H,$ implies

$x\in X_{j}=\underset{finit}{\sum }a_{i}\cdot b_{i}\subseteq I,$ we have $%
x\in I,$ so $HI\subseteq I.$

So we have $HI\subseteq H\cap I\longrightarrow (i)$

$Conversely$ suppose $x\in H\cap I$, the $x\in H$ and $x\in I.$ Since
semihyperring is multiplicatively regular, so there exist $a\in R,$ such
that, $x=x\cdot a\cdot x=x\cdot (a\cdot x)\subseteq HI,$ by using the fact
that $I$ is left hyperideal of $R$ and $(a\cdot x)\in I.$ Hence $H\cap
I\subseteq HI\longrightarrow (ii)$

From $(i)$ and $(ii)$ we have $H\cap I=HI.$

$(2)\Longrightarrow (1)$

Suppose $H\cap I=HI$ for every left hyperideal $I$ and right hyperideal $H$
of semihyperring $R.$ To prove $R$ is multiplicatively regular.

Let $a\in R,$ consider the right hyperideal $H$ of $R$ generated by $a.$
That is,

$H=aR$ and left hyperideal $I$ of $R$ generated by $a.$ That is,

$I=Ra.$

Since $1_{R}\in R,$ so obviously $a\in H$ and $a\in I.$ So $a\in H\cap I=HI$
using $(ii).$

so $a\in (aR)(Ra)=HI=\underset{b_{i}\in I}{\underset{h_{i}\in H}{\cup }}%
\left( \underset{finit}{\sum }h_{i}\cdot b_{i}\right) ,$ implies $a\in \cup
X_{i},$

where each $X_{i}=\underset{finit}{\sum }h_{i}\cdot b_{i}$ for $h_{i}$ $\in
H $ and $b_{i}\in I.$So $x\in X_{j}$ for some $j.$ So $a\in \underset{finit}{%
X_{j}=\sum }h_{i}\cdot b_{i}$ for

$h_{i}\in H=aR$ and $b_{i}\in I=Ra,$ then there exist $r_{i},$ $\acute{r}%
_{i}\in R,$ such that $h_{i}\in a\cdot r_{i},$ and

$b_{i}\in \acute{r}_{i}\cdot a,$ implies $a\in \underset{finit}{\sum }%
(a\cdot r_{i})(\acute{r}_{i}\cdot a)$ implies

$a\in (a\cdot r_{1}+a\cdot r_{2}+a\cdot r_{3}+......+a\cdot r_{n})(\acute{r}%
_{1}\cdot a+\acute{r}_{2}\cdot a+\acute{r}_{3}\cdot a+......+\acute{r}%
_{n}\cdot a)$

\qquad \qquad \qquad $a\in a\cdot (r_{1}+r_{2}+r_{3}+......+r_{n})(\acute{r}%
_{1}+\acute{r}_{2}+\acute{r}_{3}+......+\acute{r}_{n})\cdot a$

\qquad \qquad \qquad implies $a\in a\cdot (\underset{finit}{\sum }r_{i})(%
\underset{finit}{\sum }\acute{r}_{i})\cdot a$

\qquad \qquad \qquad implies $a\in a\cdot (AB)\cdot a,$ where $A=(\underset{%
finit}{\sum }r_{i})$, and $B=(\underset{finit}{\sum }\acute{r}_{i})$

\qquad \qquad \qquad implies $a\in a\cdot X\cdot a,$ where $X=AB,$ so there
exist $b\in X\subseteq R,$ such that $a\in a\cdot b\cdot a,$ implies

$a$ is multiplicatively regular element. Since $a\in R$ was arbitrary
element of $R,$ so $R$ is multiplicatively regular semihyperring.
\end{proof}

\bigskip

\section{PRIME AND SEMIPRIME HYPERIDEALS IN SEMIHYPERRING}

\subsection{Definition:}

A hyperideal $P$ of a semihyperring $R$ is called prime hyperideal of $R$,
if for hyperideals $I$ and $\ J$ of $R$ satisfying,\qquad $IJ\subseteq P,$
implies, either $I\subseteq P,$ or $J\subseteq P.$

\subsection{Proposition:}

The following conditions on hyperideal $I$ of semihyperring $R$ with
identity are equivalent.

\begin{enumerate}
\item $I$ is prime hyperideal of $R.$

\item $\{a\cdot r\cdot b:r\in R\}\subseteq I$ if and only if, either $a\in I$
or $b\in I.$

\item If $a,$ and $b$ are elements of $R,$ satisfying ;$\qquad \langle
a\rangle \langle b\rangle \subseteq I,$ then either $a\in I$ or $b\in I.$
\end{enumerate}

\begin{proof}
$(1)\Longrightarrow (2)$

Let $I$ be a prime hyperideal of $R,$ let $a,$ $b\in R,$ and the set, $%
\acute{I}=\{a\cdot r\cdot b:r\in R\}=aRb$

If $a\in I,$ or $b\in I,$ then by definition of hyperideal, $a\cdot r\cdot
b\in I,$ for all $r\in R,$ which implies that $\{a\cdot r\cdot b:r\in
R\}\subseteq I.$

$Conversely$ suppose $H=\langle a\rangle ,$ and $K=\langle b\rangle ,$ are
hyperideals of semihyperring $R$, generated by two nonzero elements of $R$,

$a,$ and $b$ respectively.

Since $1_{R}\in R$ so $\acute{I}=\{a\cdot r\cdot b\colon r\in
R\}=aRb\subseteq HK=(RaR)(RbR)$

Also since $aRb\subseteq I,$ implies $RaRbR\subseteq RIR\subseteq I$

\qquad \qquad \qquad \qquad \qquad implies $(RaR)(RbR)\subseteq I$

\qquad \qquad \qquad \qquad \qquad implies $HK\subseteq I,$

\qquad \qquad \qquad \qquad \qquad implies $H\subseteq I$, $K\subseteq I,$
implies $a\in I$ or $b\in I.$

$(2)\Longrightarrow (3)$

Let $a,$ $b\in R,$ such that $\langle a\rangle \langle b\rangle \subseteq I$

Since $\langle a\rangle =RaR,$ and $\langle b\rangle =RbR,$ so $\langle
a\rangle \langle b\rangle =(RaR)(RbR)\subseteq I$

Also since $\{a\cdot r\cdot b\colon r\in R\}\subseteq (RaR)(RbR)\subseteq I$

implies that $\{a\cdot r\cdot b\colon r\in R\}\subseteq I,$ so by $(2),$ we
have either $a\in I$ or $b\in I.$

$(3)\Longrightarrow (1)$

Let $H,$ and $K$ be the hyperideals of semihyperring $R,$ such that,

\qquad \qquad \qquad \qquad \qquad \qquad \qquad \qquad \qquad $HK\subseteq
I $

Suppose $H\nsubseteq I,$ then there exist $a\in H,$ such that $a\notin I.$
Then $\langle a\rangle \subseteq H.$

Let $b\in K,$ then $\langle b\rangle \subseteq K,$ now since $\langle
a\rangle \langle b\rangle \subseteq HK\subseteq I,$ implies $\langle
a\rangle \langle b\rangle \subseteq I.$

By $(3)$ either $a\in I$ or $b\in I.$ But $a\notin I,$ so $b\in I.$

Since $b$ was any arbitrary element of $K,$ so $K\subseteq I.$

Hence $I$ is prime hyperideal.
\end{proof}

\subsection{Corollary:}

If $a,$ and $b$ are elements of a semihyperring $R.$ Then the following
conditions on prime hyperideal $I$ of semihyperring $R$ are equivalent.

\begin{enumerate}
\item If $a\cdot b\in I,$ then $a\in I$ or $b\in I.$

\item If $a\cdot b\in I,$ then $b\cdot a\in I.$
\end{enumerate}

\begin{proof}
$(1)\Longrightarrow (2)$

Suppose $a\cdot b\in I,$ then by $(1),$ either $a\in I$ or $b\in I.$

Since $I$ is hyperideal so $b\cdot a\subseteq I.$

$(2)\Longrightarrow (1)$

Suppose $a\cdot b\in I,$ implies $(a\cdot b)\cdot r\in I,$ for all $r\in R$\
implies $a\cdot (b\cdot r)\in I,$ for all $r\in R.$

By using $(2)$ we have $(b\cdot r)\cdot a\in I,$ for all $r\in R.$implies $%
b\cdot r\cdot a\subseteq I,$for all $r\in R.$ since $I$ is prime.

So by previous Proposition, either $a\in I$ or $b\in I.$
\end{proof}

\subsection{Theorem:}

A hyperideal $I$ of a commutative semihyperring $R$ is prime if and only if $%
a\cdot b\in I,$ implies $a\in I$ or $b\in I,$ for all $a,$ $b\in R.$

\begin{proof}
Suppose $I$ is prime and $a\cdot b\subseteq I,$ implies $(a\cdot b)\cdot
r\subseteq I$ for all $r\in R,$ since $I$ is a hyperideal of $R.$ Since $R,$
is commutative

semihyperring, so we have $a\cdot r\cdot b\in I,$ for all $r\in R.$

\qquad \qquad This implies $a\cdot r\cdot b\in I.$ By Proposition 4.2 either 
$a\in I$ or $b\in I.$

$Conversely$ assume that $a\cdot b\in I,$ implies $a\in I$ or $b\in I$ for
all $a,$ $b\in R.$

Suppose $AB\subseteq I,$ where $A$ and $B$ are hyperideals of $R.$ Suppose $%
A\nsubseteq I,$ then there exists $a\in A$ such that $a\notin I.$ Now for
each $b\in B,$ $a\cdot b\in AB\subseteq I.$ By hypothesis either $a\in I$ or 
$b\in I,$ but $a\notin I,$ so $b\in I.$ Thus $B\subseteq I.$ Hence $I$ is
prime hyperideal of $R.$
\end{proof}

\subsection{Definition:}

A nonempty subset $A$ of a semihyperring $R$ is an $m-system$ if $a,$ $b\in
A,$ there exist an element $r\in R,$ such that $a\cdot r\cdot b\in A.$

\subsection{Proposition:}

A hyperideal $I$ of a semihyperring $R$ is a prime hyperideal if and only if 
$R\backslash I$ is an $m-system.$

\begin{proof}
Suppose $I$ is a prime hyperideal.

Let $a,$ $b\in R\backslash I,$ implies $a\notin I,$ and $b\notin I,$ so
Proposition$4.2$, we have$\left\{ a\cdot r\cdot b:r\in R\right\} \nsubseteq
I,$ so there exist some $r\in R,$ such that $a\cdot r\cdot b\notin I$ this
implies that $a\cdot r\cdot b\in R\backslash I,$ implies $R\backslash I$ is
an $m-system.$

$Conversely$ suppose $R\backslash I$ is an $m-system.$

Suppose $\left\{ a\cdot r\cdot b:r\in R\right\} \subseteq I.$ Suppose $%
a\notin I,$ and $b\notin I,$ implies $a\in R\backslash I,$ and $b\in
R\backslash I.$

Since $R\backslash I$ is an $m-system,$ so there exist $r_{1}\in R,$ such
that, $a\cdot r_{1}\cdot b\in R\backslash I.$

implies $a\cdot r_{1}\cdot b\notin I,$ so we have $\left\{ a\cdot r\cdot
b:r\in R\right\} \nsubseteq I,$ which is contradiction.

Hence either $a\in I$ or $b\in I,$ implies $I$ is prime hyperideal

.
\end{proof}

\subsection{Proposition:}

Any maximal hyperideal of a semihyperring $R$ with unity is prime hyperideal.

\begin{proof}
Let $M$ be a maximal hyperideal of a semihyperring $R,$ and let $\left\{
a\cdot r\cdot b:r\in R\right\} \subseteq M.$

Suppose $a\notin M,$ then $M$ is properly contained in $M+\langle a\rangle .$
But $M$ is maximal so $M+\langle a\rangle =R$

Since $1_{R}\in R,$ so we have $1_{R}\in m+s\cdot a\cdot r,$ for some $m\in
M,$ and $r,$ $s\in R,$ this implies $1_{R}\cdot b\subseteq (m+s\cdot a\cdot
r)\cdot b=m\cdot b+s\cdot (a\cdot r\cdot b)\subseteq M$

implies $b\in M,$ hence $M$ is a prime hyperideal of $R$

$.$
\end{proof}

\subsection{Proposition:}

If $I$ is hyperideal of a semihyperring $R,$ and if $H$ is hyperideal of $R,$
minimal among those ideals of $R$ properly containing $I$. Then \ $K=\{r\in
R|rH\subseteq I\}$ is prime hyperideal of $R.$

\begin{proof}
$K$ is nonempty because , for $o\in R,$ we have $o\cdot h=o,$ for all $h\in
H,$ implies $o\in K.$ So $K\neq \phi .$

Let $a,$ $b\in K,$ implies $aH\subseteq I,$ and $bH\subseteq I.$

Now $(a+b)\cdot H=a\cdot H+b\cdot H\subseteq I+I=I$

So $(a+b)\cdot H\subseteq I,$ implies $a+b\subseteq K.$

Now let $r\in R,$ and $a\in K,$ implies $aH\subseteq I$ then, $%
r(aH)\subseteq I,$ since $I$ is a hyperideal of $R.$ Also $%
r(aH)=(ra)H\subseteq I,$ implies $r\cdot a\in K.$

Now for $r\in R,$ and $a\in K,$ implies $aH\subseteq I$ then, $%
(ar)H=a(rH)\subseteq aH\subseteq I,$ because $H$ is a hyperideal and $%
rH\subseteq H,$ for all $r\in R.$ This implies $(ar)H\subseteq I,$ and hence 
$a\cdot r\in K.$ Hence $K$ is hyperideal of $R.$

Let $A,$ $B$ be hyperideals of semihyperring $R$ satisfying $AB\subseteq K.$

Assume that $B\nsubseteq K,$ to prove $K$ is prime, we must show that $%
A\subseteq K.$

Indeed since $AB\subseteq K,$ and $B\nsubseteq K,$ we have $ABH\subseteq I,$
and $BH\nsubseteq I.$ Therefore $I\subset I+BH\subseteq H,$ by minimality of 
$H,$ we have $I+BH=H,$ implies that $AI+ABH=AH,$ implies $AH=AI+ABH\subseteq
I,$ and this implies $AH\subseteq I,$ so we have $A\subseteq K.$ Hence $K$
is prime hyperideal of $R.$
\end{proof}

\subsection{Definition:}

A hyperideal $I$ of a semihyperring $R$ is called semiprime hyperideal if
for any hyperideal $H$ of $R,$ satisfying $H^{2}\subseteq I,$ implies $%
H\subseteq I.$

\subsection{Remark:}

Prime hyperideals are surely semiprime hyperideals.

\begin{proof}
Suppose $I$ is prime hyperideal of a semihyperring $R,$ then for any
hyperideals $H,$ $K,$ satisfying \ $HK\subseteq I,$ implies $H\subseteq I,$
or $K\subseteq I.$

In particular for $H^{2}\subseteq I,$ implies $H\subseteq I.$ So $I$ is
semiprime hyperideal.
\end{proof}

\subsection{Proposition:}

The following conditions on a hyperideal $I$ of a semihyperring $R$ (with
unity) are equivalent.

\begin{enumerate}
\item $I$ is semiprime .

\item $\{a\cdot r\cdot a:r\in R\}\subseteq I,$ if and only if $a\in I.$
\end{enumerate}

\begin{proof}
$(1)\Longrightarrow (2)$

Let $a\in R$ and, let $\acute{I}=\{a\cdot r\cdot a:r\in R\}$

If $a\in I,$then $a\cdot r\cdot a\in I,$ for all $r\in R,$ because $I$ is a
hyperideal.

So $\acute{I}=\{a\cdot r\cdot a:r\in R\}\subseteq I.$

$Conversely$ let $H=\langle a\rangle =RaR,$ then since $aRa\subseteq I$

implies $RaRaR\subseteq RIR\subseteq I$ implies $\left( RaR\right) \left(
RaR\right) \subseteq I$

implies $RaR\subseteq I,$ since $I$ is semiprime hyperideal, this implies $%
a\in I.$

$(2)\Rightarrow (1)$

Let $H$ be a hyperideal of $R,$ satisfying $H^{2}\subseteq I,$ let $a\in H.$
Then $\{a\cdot r\cdot a:r\in R\}\subseteq H^{2}\subseteq I,$ implies $%
\{a\cdot r\cdot a:r\in R\}\subseteq I,$ by $(2),$ we have $a\in I,$ since $a$
was arbitrary element of $H,$ so $H\subseteq I,$ and we have $I$ is
semiprime hyperideal.
\end{proof}

\subsection{Definition:}

A nonempty subset $A$ of a semihyperring $R$ is a $p-system$ if and only if
for $a\in A,$ there exists an element $r\in R,$ such that $a\cdot r\cdot
a\in A.$

\subsection{Corollary:}

A hyperideal $I$ of a semihyperring $R$ is semiprime hyperideal if and only
if $R\backslash I$ is $p-system.$

\begin{proof}
Suppose $I$ is semiprime hyperideal of a semihyperring $R.$

Let $a\in R\backslash I,$ implies $a\notin I,$ then by above Proposition we
have $\{a\cdot r\cdot a:r\in R\}\nsubseteq I,$ so there exist $r_{1}\in R,$
such that $a\cdot r_{1}\cdot a\notin I\ $that is $a\cdot r_{1}\cdot a\in
R\backslash I,$ implies $R\backslash I$ is a $p-system.$

$Conversely$ suppose $R\backslash I$ is a $p-system.$

\qquad \qquad \qquad Let $\{a\cdot r\cdot a:r\in R\}\subseteq I.$

To prove $a\in I,$ on contrary suppose that $a\notin I,$ implies $a\in
R\backslash I.$ Since $R\backslash I$ is a $p-system,$ so there exists an
element $r_{1}\in R,$ such that $a\cdot r_{1}\cdot a\in R\backslash I,$
implies $a\cdot r_{1}\cdot a\notin I,$ and this implies that $\{a\cdot
r\cdot a:r\in R\}\nsubseteq I,$ which is a contradiction. So we have $a\in
I. $ Hence $I$ is a semiprime hyperideal.
\end{proof}

\subsection{Remark:}

Any $m-system$ is a $p-system.$

\begin{proof}
Suppose $A$ is an $m-system$ of a semihyperring $R.$ Then for $a,$ $b\in A,$
there exist $r\in R,$ such that $a\cdot r\cdot b\in A,$ in particular for $%
a\in A,$ $a\cdot r\cdot a\in A.$

Hence $A$ is a $p-system.$
\end{proof}

\bigskip

\section{IRREDUCIBLE AND STRONGLY IRREDUCIBLE \ HYPERIDEALS}

\subsection{Definition:}

A hyperideal $I$ of a semihyperring $R$ is irreducible if and only if for
ideals $H$ and $K$ of $R$ $,$ $I=H\cap K,$ implies $I=H,$ or $I=K.$ The
hyperideal $I$ is strongly irreducible if and only if for ideals $H$ and $K$
of $R,$ $H\cap K\subseteq I$ implies $H\subseteq I$ or $K\subseteq I.$

\subsection{Remark:}

Every strongly irreducible hyperideal is an irreducible hyperideal.

\begin{proof}
Let $I$ be a strongly irreducible hyperideal of $R$. Let $I=H\cap K\subseteq
I,$ implies $H\subseteq I$ or $K\subseteq I.$ Also $I=H\cap K\subseteq H,$
and $I=H\cap K\subseteq K,$ implies $I\subseteq H,$ and $I\subseteq K.$ So
we have $I=H,$ or $I=K.$
\end{proof}

\subsection{Definition:}

A non-empty set $A$ of a semihyperring $R$ is an $i-system$ if and only if
for $a,$ $b\in A,$ implies $\langle a\rangle \cap \langle b\rangle \cap
A\neq \phi .$

\subsection{Proposition:}

The following conditions on a hyperideal $I$ of a semihyperring $R$ are
equivalent,

\begin{enumerate}
\item $I$ is strongly irreducible.

\item If $a,$ $b\in R,$ satisfying $\langle a\rangle \langle b\rangle
\subseteq I,$ implies $a\in I,$ or $b\in I.$

\item $R/I$ is an $i-system.$
\end{enumerate}

\begin{proof}
$(1)\rightarrow (2)$

Let for $a,$ $b\in R,$ we have $\langle a\rangle \langle b\rangle \subseteq
I,$ since $I$ is strongly irreducible so $\langle a\rangle \subseteq I,$ or $%
\langle b\rangle \subseteq I.$

Which implies that $a\in I,$ or $b\in I.$

$(2)\rightarrow (3)$

Suppose if $a,$ $b\in R,$ satisfying $\langle a\rangle \langle b\rangle
\subseteq I,$ implies $a\in I,$ or $b\in I.$ To prove $R/I$ is an $i-system.$

On contrary suppose $R/I$ is not an $i-system,$ then for $a,$ $b\in R/I$,
implies $(\langle a\rangle \cap \langle b\rangle )\cap R/I=\phi ,$ so we get
from here that $\langle a\rangle \langle b\rangle \subseteq I,$ this implies
that $a\in I,$ or $b\in I,$ which is contradiction to the fact that $a,$ $%
b\in R/I.$ Thus we have $(\langle a\rangle \cap \langle b\rangle )\cap
R/I\neq \phi .$ So $R/I$ is an $i-system.$

$(3)\rightarrow (1)$

Suppose $R/I$ is an $i-system.$

To prove $I$ is strongly irreducible.

Let $H$ and $K$ be the ideals of $R$ satisfying $H\cap K\subseteq I.$ To
prove $H\subseteq I$ or $K\subseteq I.$

Suppose neither $H\subseteq I$ nor $K\subseteq I,$ then there exists $a,$ $%
b\in R$ such that $a\in H\backslash I$ and $b\in K\backslash I,$ which
implies $a,$ $b\in R\backslash I,$ also since $R\backslash I$ is an $%
i-system $ so $\langle a\rangle \cap \langle b\rangle \cap R/I\neq \phi ,$
so there exists $c\in \langle a\rangle \cap \langle b\rangle \cap R/I,$
implies $c\in R\backslash I.$ But $c\in \langle a\rangle \cap \langle
b\rangle \subseteq H\cap K,$ implies $H\cap K\nsubseteq I.$

Which is a contradiction. Hence either $H\subseteq I$ or $K\subseteq I.$
\end{proof}

\subsection{Remark:}

Any prime hyperideal is strongly irreducible.

\begin{proof}
Let $I$ be a prime hyperideal of a semihyperring $R.$

Let $H$ and $K$ be ideals of semihyperring $R,$ such that, $H\cap K\subseteq
I.$

Since $HK\subseteq H\cap K\subseteq I,$ implies $HK\subseteq I,$ and $I$ is
prime so either $H\subseteq I$ or $K\subseteq I.$ Hence for $H\cap
K\subseteq I,$ we have either $H\subseteq I$ or $K\subseteq I.$

$I$ is strongly irreducible.
\end{proof}

\subsection{Proposition:}

Let $0\neq a\in R,$ and $I$ be a hyperideal of $R$ not containing $a$, then
there exists an irreducible hyperideal $H$ of $R,$ which contains $I,$ but
not contain $I.$

\begin{proof}
Let $\mathit{A}$ be the collection of all ideals of $R$ which contains $I$,
but not contain $a,$ that is;

\qquad \qquad $\mathit{A}=\{hyperideal$ $J$ $%
{\vert}%
$ $I\subseteq J,$ $a\notin J\}$

Clearly $\mathit{A}\neq \Phi ,$ because $I\subseteq I,$ and $a\notin I.$

Let $\{H_{i}$ $%
{\vert}%
$ $i\in \wedge \}$ be a chain in $\mathit{A,}$ and consider $\acute{H}=%
\underset{i\in \wedge }{\cup }H_{i}.$

Let $a,$ $b\in \acute{H},$ then $a\in H_{i},$ and $b\in H_{j}$ for some $i$
and $j$ in $\wedge .$ Then either $H_{i}\subseteq H_{j}$ or $H_{j}\subseteq
H_{i}.$ Without any loss of generality

we assume that $H_{i}\subseteq H_{j},$ then $a,$ $b\in H_{j},$ implies $%
a+b\subseteq H_{j}\subseteq \acute{H},$ and this implies $a+b\subseteq 
\acute{H}.$

Now also for $r\in R,$ and for any $a\in \acute{H},$ implies $a\in H_{i},$
for some $i$ in $\wedge ,$ then $r\cdot a,$ $a\cdot r\in H_{i}\subseteq 
\acute{H}.$

Hence $\acute{H}$ is hyperideal of $R.$ Since each $Hi\in \mathit{A}$
contains $I$ but does not contain $a,$ so $H$ also contains $I$ but does not
contain $a.$ Thus $\acute{H}$ is an upper bound for the chain $\{H_{i}$ $%
{\vert}%
$ $i\in \wedge \}.$ Therefore by $Zorn^{\prime }s$ $lemma$ $\mathit{A}$
contains a maximal element say $H.$

Now we show that $H$ is an irreducible hyperideal. Let $J$ and $K$ be the
hyperideals of $R,$ with \ $J\cap K=H.$

Suppose $J\neq K$ and $K\neq H.$ As $H\subseteq J$ and $\ H\subseteq K,$ so
by the maximality of $H,$ $a\in K$ and $a\in J.$

Thus $a\in J\cap K=H,$ implies $a\in H,$ which is a contradiction. Hence
either $H=K$ or $H=J.$ So $H$ is irreducible.
\end{proof}

\subsection{Proposition:}

Any hyperideal of a semihyperring $R$ is the intersection of all irreducible
hyperideals of $R$ containing it.

\begin{proof}
Suppose $\acute{I}$ be the intersection of all irreducible ideals of $R$
which contains $I,$ that is \ $I\subseteq \acute{I}.$

Let $0\neq a\in R,$ such that $a\notin I,$ then by above Proposition, there
exist an irreducible hyperideal $H$ of $R$ containing $I,$ but not
containing $a$ that is $I\subseteq H$ and $a\notin H.$ $H$ is present in the
collection, whose intersection is $\acute{I}$, so $a\notin \acute{I},$
implies $\acute{I}\subseteq I.$ So we have $I=\acute{I}.$
\end{proof}

\subsection{Proposition:}

A hyperideal $I$ of a semihyperring $R$ is prime if and only if it is
semiprime and strongly irreducible.

\begin{proof}
Suppose $I$ is prime hyperideal of a semihyperring$R.$ Then obviously $I$ is
semiprime hyperideal.

Suppose $H$ and $K$ are hyperideal of $R$, such that \qquad $H\cap
K\subseteq I.$

Since $HK\subseteq H\cap K\subseteq I,$ implies $HK\subseteq I.$ Since $I$
is prime hyperideal\ so either $H\subseteq I,$ or $K\subseteq I.$

Hence $I$ is strongly irreducible.

$Conversely$ assume that $I$ is semiprime and strongly irreducible
hyperideal.

Let $A$ and $B$ be hyperideals of semihyperring $R,$ such that $AB\subseteq
I.$

Since $A\cap B\subseteq A$ and $A\cap B\subseteq B,$ implies $(A\cap
B)^{2}\subseteq AB\subseteq I,$ implies $(A\cap B)^{2}\subseteq I.$ Since $I$
is semiprime so $A\cap B\subseteq I.$ Also $I$ is strongly irreducible, so
we have $A\subseteq I$ or $B\subseteq I.$ Hence $I$ is prime.
\end{proof}

\subsection{Definition:}

We denote by $I(R)$, the lattice of hyperideals of a semihyperring $R$ and
by $H(R)$ the set of irreducible hyperideals of $R$. For any hyperideal $I$
of $R$, we define

$\Theta _{I}=\{J\in H(R):I\varsubsetneq J\}$

and $\tau (H(R))=\{\Theta _{I}:I\in I(R)\}$

\subsection{Theorem:}

The set $\tau (H(R))$ form a topology on the set $H(R).$

\begin{proof}
\bigskip First we show that $\tau (H(R))$ form a topology on a set H(R).

As $\Theta _{<o>}=\{J\in H(R):<o>\varsubsetneq J\}=\varphi ,$ since $<o>$ is
contained in every (irreducible) hyperideal of $R$.

So $\Theta _{<o>}$ is an empty subset of $\tau (H(R))$. On the other hand $%
\Theta _{R}=\{J\in H(R):R\varsubsetneq J\}=H(R)$ is also true since

irreducible hyperideals are proper. So $\Theta _{R}=H(R)$ is also an element
of $\tau (H(R)).$

Now let $\Theta _{I}\cap \Theta _{M}=\{J\in H(R):I\varsubsetneq J$ and $%
M\varsubsetneq J\}=\{J\in H(R):I\cap M\varsubsetneq J\}=\Theta _{I\cap M}.$

Hence $\Theta _{I}\cap \Theta _{M}$ is also an element of $\tau (H(R)).$

Let us consider an arbitrary family $\{I_{\lambda }:\lambda \in \Lambda \}$
of hyperideals of $R$.

Since $\cup \Theta _{I_{\lambda }}=\cup \{J\in H(R):I_{\lambda
}\varsubsetneq J\}=\{J\in H(R):\exists \lambda \in \Lambda $ so that $%
I_{\lambda }\varsubsetneq J\}=\{J\in H(R):\sum I_{\lambda }\varsubsetneq
J\}=\Theta _{\sum I_{\lambda }}.$

Since $_{\sum I_{\lambda }}\in I(R)$, it follows that $\cup \Theta
_{I_{\lambda }}\in \tau (H(R)).$ Thus the set $\tau (H(R))$ form a topology
on $H(R).$
\end{proof}

\subsection{Remark:}

The assignment $I\longrightarrow \Theta _{I}$ is an isomorphism between the
lattice $I(R)$ of hyperideals of semihyperring $R$ and the lattice of open
subsets of $H(R)$.

\subsection{Definition:}

\ The set $H(R)$ of irreducible hyperideals of $R$ is called irreducible
spectrum of $R$. The topology $\tau (H(R))$ in above theorem is called the
irreducible spectrum topology on $H(R)$. $H(R)$ is called irreducible
spectral space.

\subsection{Proposition:}

Let $R$ be a commutative semihyperring, then the following are equivalent.

\begin{enumerate}
\item $R$ is multiplicatively regular.

\item Every hyperideal of $R$ is idempotent.

\item $I\cap J=IJ$ for every hyperideal $I$ and $J$ of $R.$

\item Every hyperideal of $R$ \ is semiprime.
\end{enumerate}

\begin{proof}
$(1)\Longrightarrow (2)$

Let $R$ be multiplicatively regular semihyperring.

Let $I$ be any hyperideal of $R.$ Since $I^{2}=I\cdot I\subseteq I\cdot
R\subseteq I,$ implies $I^{2}\subseteq I.$ Let $x\in I$, since $R$ is
regular semihyperring, so there exists $a\in R,$ such that \qquad $x=xax.$

Since $x\in I,$ implies $a\cdot x\in I,$ for all $a\in R.$ Also since $%
x=xax=x(ax)\in I\cdot I=I^{2},$ implies $x\in I^{2}.$ Hence $I\subseteq
I^{2},$ so we have finally $I=I^{2}.$ Hence every hyperideal of $R$ is
idempotent.

$(2)\Longrightarrow (3)$

Suppose every hyperideal of $R$ is idempotent.

Let $I$ and $J$ be the hyperideals of $R.$ Then obviously $IJ\subseteq I\cap
J.$

But also $I\cap J\subseteq I,$ and $I\cap J\subseteq J.$ Since $I\cap J$
being intersection of two hyperideals is again an hyperideal. So $I\cap
J=(I\cap J)^{2}\subseteq IJ.$

Hence we have $I\cap J=IJ.$

$(3)\Longrightarrow (4)$

Suppose $I\cap J=IJ$ for every hyperideal $I$ and $J$ of semihyperring $R.$

Let $I$ and $J$ be any hyperideals of $R,$ such that ,\qquad $I^{2}\subseteq
J.$

\qquad Since $I=I\cap I=I\cdot I\subseteq J,$ implies $I\subseteq J.$ Hence $%
J$ is semiprime hyperideal.

$(4)\Longrightarrow (3)$

Suppose every hyperideal of semihyperring is semiprime.

To prove $I\cap J=IJ$ for every hyperideal $I$ and $J$ of semihyperring $R.$

Let $I$ and $J$ be the hyperideals of $R.$ Then obviously $IJ\subseteq I\cap
J.$

But also $I\cap J\subseteq I,$ and $I\cap J\subseteq J.$ Since $I\cap J$
being intersection of two hyperideals is again a hyperideal, and every

hyperideal of semihyperring is semiprime. So we have $(I\cap J)^{2}\subseteq
IJ,$ implies $(I\cap J)\subseteq IJ.$

Hence $I\cap J=IJ.$

$(3)\Longrightarrow (1)$

Suppose $H\cap I=HI$ for every left hyperideal $I$ and right hyperideal $H$
of semihyperring $R.$ To prove $R$ is multiplicatively regular.

Let $a\in R,$ consider the right hyperideal $H$ of $R$ generated by $a.$
That is,

$H=aR,$ and left hyperideal $I$ of $R$ generated by $a.$ That is,

$I=Ra.$

Since $1_{R}\in R,$ so obviously $a\in H$ and $a\in I.$ So $a\in H\cap I=HI$
using $(ii).$

so $a\in (aR)(Ra)=HI=\underset{b_{i}\in I}{\underset{h_{i}\in H}{\cup }}%
\left\{ \underset{finit}{\sum }h_{i}\cdot b_{i}\right\} ,$ implies $a\in
\cup X_{i},$

where each $X_{i}=\underset{finit}{\sum }h_{i}\cdot b_{i}$ for $h_{i}$ $\in
H $ and $b_{i}\in I.$So $x\in X_{j}$ for some $j.$ So $a\in \underset{finit}{%
X_{j}=\sum }h_{i}\cdot b_{i}$ for $h_{i}\in H=aR$ and

$b_{i}\in I=Ra,$ then there exist $r_{i},$ $\acute{r}_{i}\in R,$ such that $%
h_{i}\in a\cdot r_{i},$ and $b_{i}\in \acute{r}_{i}\cdot a,$ implies $a\in 
\underset{finit}{\sum }(a\cdot r_{i})(\acute{r}_{i}\cdot a)$ implies

\ \ \ \ \ \ \ \ \ \ \ \ \ \qquad implies $a\in a\cdot (\underset{finit}{\sum 
}r_{i})(\underset{finit}{\sum }\acute{r}_{i})\cdot a$

\qquad \qquad \qquad implies $a\in a\cdot (A\cdot B)\cdot a,$ where $A=(%
\underset{finit}{\sum }r_{i})$, and $B=(\underset{finit}{\sum }\acute{r}%
_{i}) $

\qquad \qquad \qquad implies $a\in a\cdot X\cdot a,$ where $X=A\cdot B,$ so
there exist $b\in X\subseteq R,$ such that $a=a\cdot b\cdot a,$ implies $a$
is

multiplicatively regular element. Since $a\in R$ was arbitrary element of $%
R, $ so $R$ is multiplicatively regular semihyperring.
\end{proof}

\end{document}